\documentclass[11pt,a4paper,leqno]{article}
\usepackage{graphicx}
\usepackage{ascmac}
\usepackage{hyperref}
\usepackage{amsmath,amssymb}
\usepackage{amsthm}
\usepackage{newpxtext,newpxmath}
\title {Other special cases of a square problem. }
\author {Yasushi Ieno}
\date{}
\begin{document}
\maketitle 
\noindent Abstract. There exists "a square problem": in a unit square is there a point with four rational distances to the vertices? This problem is still regarded as unproved. Yang showed proofs for several special cases of the square problem. By the reference of Yang's researches, We have proved other special cases of this problem.\\

\\\\
0. Introduction\\

45 years ago, C.W.Dodge asked "a square problem" in the Mathematical Magazine [1]: is there a point in a unit square which has all rational distances to the four vertices? Since then a lot of researchers have tried to prove this square problem. But still now this is regarded as unproved [2].

Yang proved several special cases of this problem [3]; if the point is located on diagonals, midlines, or edge, or the side length of the square is n times the distance from the point to one side, where n and $\rm n^2$+4 are both prime numbers, this point cannot have four rational distances to the vertices of the unit square.

By the reference of Yang's researches, we have proved other specail cases of this problem. For example, assuming a square is with side-width of an integer, the distance from the point to one side cannot be a prime. 
\\\\
1. Preliminary knowledge.\\

It is evident that this square problem is equivalently regarded as a problem below.\\  
\begin{shadebox}
{\bf Is there a point in a square with side-width of on integer which has all integer distances to the four vertices?}
\end{shadebox}
\quad\par
This new problem is gained by multiplying with a suitable integer. 

Then we obtain a square with side-width of integer number z, then its four vertices are (0,0), (0,z), (z,z) and (z,0) respectively. The point P(x,y) has four integer distances to the sides, and of course, four integer distances to the vertices, and both x and y are non-negative integer. Holding such properties we can assume that integer x, y and z have no common divisor more than 1. 
\\\\
2. Our researches.\\

As mentioned above, we will research a square with four vertices  (0,0), (0,z), (z,z) and (z,0), and the point P(x,y), where x, y and z are non-negative integer and they have no common divisor more than 1.

We easily see that if neither x nor y is zero, then z and one of x and y are even. Assume it is y, then x must be odd, which results in that P is not located on the diagonals. In addition, y is a multiple of 4.

Furthermore Yang showed 3 proofs as the following lemmas [3].
\\\\
\noindent {\bf Lemma 1.} 

On the edge of a square with side-width of an integer, no point has four integer distances to the vertices.
\\\\
\noindent {\bf Lemma 2.} 

On the midlines of a square with side-width of an integer, no point has four integer distances to the vertices.
\\\\
\noindent {\bf Lemma 3.} 

Considering a square with side-width of an integer, if the side-length is n times the distance from a point to one side of the square, where n and $\rm n^2$+4 are both prime numbers, no point has four integer distances to the vertices.
\\

From now on, we assume that x is odd and y is even, and naturally z is even, for because of Lemma 1, x and y are positive.

Then x and y can be expressed as $\rm x{=}k(u^2{-}v^2)$ and $\rm y{=}2kuv$ where u, v and k is positive integers. One of u and v is even, and the other is odd, otherwise x and y must be even, which contradicts the assumption. Therefore z is a multiple of 4.

\begin{figure}[htbp]
\begin{center}
\includegraphics [width=40mm]{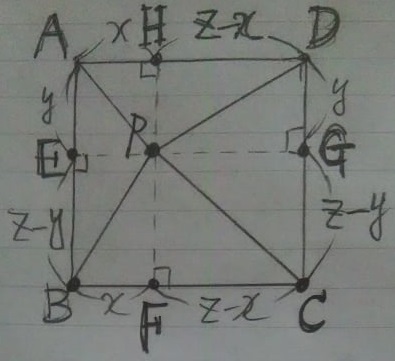}
\caption{a square}
\end{center}
\end{figure}

Now on Figure 1, $\rm \overline{AP}^2$=$\rm \overline{AH}^2$+$\rm \overline{AE}^2$=$\rm x^2$+$\rm y^2$. If x$\not\equiv$0 (mod 3) and y$\not\equiv$0 (mod 3), then $\rm \overline{AP}^2$=1+1=2, which is impossible, so one of x and y is a multiple of 3. This theory holds not only for the combination (x and y) but also for the three combinations (x and z${-}$y), (z${-}$x and z${-}$y) and (z${-}$x and y). So z is a multiple of 3, naturally a multiple of 12. \\

We check again that x is a positive odd integer, y is a positive multiple of 4 and z is a positive multiple of 12. \\

Next we are showing several theorems as follows.
\\\\
\noindent {\bf Theorem 1.} 

The following inequalities hold.

\ $\rm x^2{\geq}2y{+}1$,  $\rm y^2{\geq}2x{+}1$,  $\rm (z{-}x)^2{\geq}2y{+}1$,  $\rm y^2{\geq}2(z{-}x){+}1$, 

\ $\rm x^2{\geq}2(z{-}y){+}1$,  $\rm (z{-}y)^2{\geq}2x{+}1$,  $\rm (z{-}x)^2{\geq}2(z{-}y){+}1$ and $\rm (z{-}y)^2{\geq}2(z{-}x){+}1$.

\begin{proof}
\quad\par
On Figure 1, $\rm AP^2$=$\rm AH^2$+$\rm AE^2$=$\rm x^2$+$\rm y^2$. Now AP is as long as a positive integer, so $\rm AP{\geq}y{+}1$.\\
$\rm x^2$+$\rm y^2{\geq}{(y{+}1)}^2$, which results in that $\rm x^2{\geq}2y{+}1$. But also $\rm AP{\geq}x{+}1$, so $\rm y^2{\geq}2x{+}1$.\\
Similarly considering the case of BP, CP and DP, we easily see that the remaining six inequalities hold.

This completes Theorem 1.\\
\end{proof}

\noindent {\bf Theorem 2.} 

For an odd prime number p, if $\rm (\tfrac{2}{p}){=}{-}1$, namely 2 is not a quadratic residue modulus p, then x$\not\equiv$y(mod p) and $\rm z{-}x{\not\equiv}$y(mod p) and $\rm z{-}x{\not\equiv} z{-}y$(mod p) and $\rm z{-}x{\not\equiv}z{-}y$(mod p).

\begin{proof}
\quad\par
Assume x$\equiv$y(mod p) and $\rm (\tfrac{2}{p}){=}{-}1$. It holds that x$\equiv$y$\equiv$n(mod p) where n is a non-negative integer less than p.

If n$\neq$0, then $\rm x^2$+$\rm y^2$$\equiv$2$\rm n^2{\not\equiv}0$(mod p). Then there does not exist an integer w such that $\rm x^2$+$\rm y^2$$\equiv$$\rm w^2$(mod p), which is contradiction. So n must be equal to 0. 

This theory holds not only for the combination (x and y) but also for the three combinations (x and z${-}$y), (z${-}$x and z${-}$y) and (z${-}$x and y). 
If for one of the four combinations n=0, then for the other three combinations n must be 0, which is contradiction because x, y and z have no common divisor more than 1.

This completes Theorem 2.
%
\\
\end{proof}

\noindent {\bf Theorem 3.} 

Neither x nor $\rm z{-}x$ is an odd prime number. 
\begin{proof}
\quad\par
x and y constitute two non-hypotenuse sides of a right triangle whose side lengths are integer.

So x and y can be expressed as x=k($\rm u^2{-}v^2$)=k$\rm (u{+}v)(u{-}v)$ and y=2kuv where u, v and k are positive integer numbers, and u and v are coprime. Assume that x is an odd prime number. Then k=1 and u$\rm {-}$v=1, so x=$\rm u{+}v$=$\rm 2v{+}1$ and y=2uv=2$\rm (v{+}1)v$, which results in that y=$\rm (x^2{-}1)/2$. Now y is uniquely decided by x.

Therefore also $\rm z{-}y$ must be $\rm (x^2{-}1)/2$, which contradicts Lemma 2 because y=$\rm z{-}y$. Similarly assume that $\rm z{-}x$ is an odd prime number, which is contradiction, too.

This completes Theorem 3. 
\\
\end{proof}

\noindent {\bf Theorem 4.} 

For an odd prime number p, if $\rm (\tfrac{2}{p}){=}{-}1$, then $\rm x{\neq} p^n$ and $\rm z{-}x{\neq} p^n$, where n is a positive integer.
\begin{proof}
\quad\par
Assume $\rm x{\neq} p^n$. As mentioned above, x and y can be expressed as x=k($\rm u^2{-}v^2$)=k$\rm (u{+}v)(u{-}v)$ and y=2kuv where u, v and k are positive integer numbers, and u and v are coprime. If k is more than 1, k must be an exponential of p, now which is contradiction because of Theorem 2. So k =1. Because u and v are coprime and $\rm x{\neq} p^n$, $\rm u{-}v$=1. So like the proof of Theorem 2, y=$\rm (x^2{-}1)/2$ and y is uniquely decided by x, also $\rm z{-}y$ must be $\rm (x^2{-}1)/2$, which contradicts Lemma 2 because y=$\rm z{-}y$. Similarly assume that $\rm z{-}x$ is an odd prime number, which is contradiction, too.

This completes Theorem 4.
\\
\end{proof}

\noindent {\bf Theorem 5.} 

For odd prime numbers $\rm p_1,p_2,. . .,p_j$, if any of them is either equal to 1(mod 4) or holds that $\rm (\tfrac{2}{p_i}){=}{-}1$ for any $\rm 1{\leq}i{\leq}j$, then neither y nor $\rm z{-}y $ can be equal to $\rm 2^{h{+}1}p_1^{l_1}p_2^{l_2}. . .p_j^{l_j}$ such that $\rm 2^h{>}p_1^{l_1}p_2^{l_2},. . .,p_j^{l_j}$, where h,$\rm l_1,l_2,. . .,l_j$ are positive integers.
\begin{proof}
\quad\par
Assume $\rm y{=}2^{h{+}1}p_1^{l_1}p_2^{l_2},. . .,p_j^{l_j}$. As mentioned above, x and y can be expressed as x=k($\rm u^2{-}v^2$)\\=k$\rm (u{+}v)(u{-}v)$ and y=2kuv where u, v and k are positive integer numbers, and u and v are coprime.

At first, according to Theorem 1, we see that k cannot include as its factor any prime p, if modulus which 2 is not a quadratic residue.

So now any prime factor of k is such that 2 is a quadratic residue modulus it, and it is equal to 1 (mod 4).   
%

As the assumption above, kuv=$\rm 2^hp_1^{l_1}p_2^{l_2}. . .p_j^{l_j}$.

One of u and v is even, and the other is odd. Now u involves $\rm 2^h$ as its factor, and v does not involve 2, otherwise u$\rm \leq{p_1^{l_1}p_2^{l_2}. . .p_j^{l_j}}$ and v$\rm {\geq}2^h$, so $\rm u^2{-}v^2{<}0$, which contradicts the assumption. Therefore $\rm u^2{-}v^2{\equiv}3$ (mod 4), and k$\rm {\equiv}1$ (mod 4) because any factor of k is equal to 1 (mod 4), which results in that x=$\rm k(u^2{-}v^2){\equiv}3$ (mod 4). Similarly $\rm z{-}x{\equiv}3$ (mod 4), then $\rm z{=}x{+}(z{-}x) {\equiv}2$ (mod 4), which is contradiction, for z is a multiple of 12.

This completes Theorem 5.
\\
\end{proof}
\newpage
Then two corollaries follow.\\
 
\noindent {\bf Corollary 5.1.} 

Neither y nor $\rm z{-}y $ can be equal to $\rm 2^{h{+}1}$ where h is a positive integer.\\

\noindent {\bf Remark 5.1.}   This Corollary 5.1 is a special case for Theorem 5.\\

\noindent {\bf Corollary 5.2.} 

Neither x nor $\rm z{-}x $ can be equal to $\rm q_1q_2$ \\
such that $\rm both\ q_1\ and\ q_2\ are\ odd\ positive\ integers$, $\rm (\tfrac{2}{q_1}){=}(\tfrac{2}{q_2}){=}{-}1$,\\
and $\rm (q_1^2{-}q_2^2)/4{=}2^{h}p_1^{l_1}p_2^{l_2}. . .p_j^{l_j}$ where h,$\rm l_1,l_2,. . .,l_j$ are positive integers, $\rm 2^h{>}p_1^{l_1}p_2^{l_2},. . .,p_j^{l_j}$ and $\rm p_1,p_2,. . .,p_j$ are odd prime numbers, any of which is either equal to 1(mod 4) or holds that $\rm (\tfrac{2}{p_i}){=}{-}1$ for any $\rm 1{\leq}i{\leq}j$.\\

\noindent {\bf Corollary 5.3.} 

Neither x nor $\rm z{-}x $ can be equal to $\rm q_1q_2$ such that $\rm both\ q_1\ and\ q_2\ are\ odd\ positive\ integers$, $\rm (\tfrac{2}{q_1}){=}(\tfrac{2}{q_2}){=}{-}1$, and $\rm (q_1^2{-}q_2^2)/4{=}2^{h}$ where h is a positive integer.\\

\noindent {\bf Remark 5.3.}   This Corollary 5.3 is a special case for Corollary 5.2.\\

\noindent {\bf Theorem 6.} 

 x and $\rm z{-}x$ cannot consist of the combination of $\rm p_1p_2\ and\ q_1q_2,\ where\ p_1,p_2,q_1\ and\ q_2$ are odd prime numbers such that $\rm p_1{\neq}p_2$  , $\rm q_1{\neq}q_2$, $\rm (\tfrac{2}{p_1}){=}{-}1$, $\rm (\tfrac{2}{p_2}){=}{-}1$, $\rm (\tfrac{2}{q_1}){=}{-}1$ and $\rm (\tfrac{2}{q_2}){=}{-}1$.
\begin{proof}
\quad\par
Assume $\rm x{=}p_1p_2$ and $\rm z{-}x{=}q_1q_2$ where $\rm p_1{>}p_2$ and  $\rm q_1{>}q_2$. As mentioned above, \\
x and y can be expressed as x=k($\rm u^2{-}v^2$)=k$\rm (u{+}v)(u{-}v)$ and y=2kuv where u, v and k are positive integer numbers, and u and v are coprime. Now $\rm (\tfrac{2}{p_1}){=}{-}1$ and $\rm (\tfrac{2}{p_2}){=}{-}1$, so k=1. Then $\rm x{=}p_1p_2{=}(u{+}v)(u{-}v)$, which results in that  (u,v)= 
(($\rm p_1{+}p_2)/2,(p_1{-}p_2)/2$) or (($\rm p_1p_2{+}1)/2,(p_1p_2{-}$1)/2). And now y=2uv, so y and $\rm z{-}y$ consist of the combination of $\rm (p_1^2{-}p_2^2)/2\ and\ (p_1^2p_2^2{-}1)/2$. 

Similar for the case of $\rm z{-}x$ and y, we see that y and $\rm z{-}y$ consist of the combination of $\rm (q_1^2{-}q_2^2)/2\ and\ (q_1^2q_2^2{-}1)/2$.

Therefore $\rm (p_1^2{-}p_2^2)/2{=}(q_1^2{-}q_2^2)/2$ and $\rm (p_1^2p_2^2{-}1)/2{=}\rm (q_1^2q_2^2{-}1)/2$, which results in $\rm p_1{=}q_1$ and $\rm p_2{=}q_2$. We see that $\rm x{=}z{-}x$, which is contradiction, owing to Lemma 2.

This completes Theorem 6.
\\
\end{proof}
\quad\par
\noindent 3. Discussion and Conclusions.\\

Yang proved that there are some points, expressed as a relative position of a square, which cannot be P(x,y) fulfilling the given condition of the square problem. 

In contrast, we have proved that there are some points, expressed absolutely as a concrete position with the values of x and y, which similarly cannot be P(x,y) fulfilling the given condition of the square problem, on the basis of Yang's paper.

As mentioned at the beginning, this square problem has not yet been perfectly proved, but there have been negative views regarding the existence of such points, though it is known that there are innumerable points inside a unit square such that as many as three distances to the three vertices are rational numbers [4], for example the cases  (x,y,z)=(7,24,52), (297,304,700) and so on.

And in this paper unwillingly we also have not proved this problem perfectly yet.

But we can easily see that many integers are unavailable for x or y, under the condition that  x must be odd and that y must be a multiple of 4.\\

For example, \\
when z=60, according to the theorems above, we see as follows. 

By Theorem 3,\\
for x, integers such as 3, 5, 7, 11, 13, 17, 23, 29, 31, 37, 41, 47, 53 and 59 are unavailable, and the same for $\rm z{-}x$, \\
by symmetry for x, integers such as 1, 3, 5, 7, 11, 13, 17, 19, 23, 29, 31, 37, 41, 43, 47, 49, 53, 55, 57 and 59 are unavailable, consequently.

By Theorem 4,\\
for x, integers such as 3, 5, 9, 25 and 27 are unavailable, and the same for $\rm z{-}x$, \\
by symmetry for x, integers such as 3, 5, 9, 25, 27, 33, 35, 51, 55 and 57 are unavailable, consequently.

By Theorem 5,\\
for y, integers such as 4, 8, 16, 24, 32 and 48 are unavailable, and the same for $\rm z{-}y$, \\
by symmetry for y, integers such as 4, 8, 12, 16, 24, 28, 32, 36, 44, 48, 52 and 56 are unavailable, consequently; 20 and 40 are still decided to be unavailable for y. But according to Lemma 3, for the case y=20 and y=40 then the side-length is 60, so 3 times the distance from a point to one side of the square, 20. And both 3 and $3^2{+}4$=13 are prime. Therefore 20 and 40 are unavailable. As a result, we see that for the case z=60 there is no available pair of x and y for the problem, in fact, without checking Theorem 3 through Theorem 4.\\

As mentioned above, the theorems we have shown are effective tools for a concrete case to check rapidly that there does not exist an available point P or not. But they are far from a perfect proof. \\

We will continue to make efforts to approach the proof of this square problem.
\\\\
\newpage
\centerline{References}
\quad\par
\noindent [1] C. W. Dodge, Math. Mag. 49, 43 (1976).\\
\noindent [2] Peter Brass, William Moser, Janos Pach, Research Problems in Discrete Geometry, Springer, 2005.\\
\noindent [3] Yang Ji, arvix, available at \href{https://arxiv.org/abs/2105.05250/}{https://arxiv.org/abs/2105.05250}\\
\noindent [4] A215365 of the OEIS, available at \href{https://oeis.org/A215365}{https://oeis.org/A215365}\\
\end{document}